\documentclass[11pt]{article}

\usepackage{fullwidth}
\usepackage{xypic}
\usepackage[top=1.5in, bottom=1.5in, left=1in, right=1in]{geometry}
\usepackage{amsgen}
\usepackage{amsmath}
\usepackage{amstext}
\usepackage{amsbsy}
\usepackage{amsopn}
\usepackage{amsfonts}
\usepackage{amssymb}
\usepackage{eepic}
\usepackage{graphicx}
\usepackage{epsf}
\usepackage{pstricks}
\xyoption{all}

\def\Box{\square}
\def\edge{\relbar\joinrel\relbar}

\def\tra#1{\smash{\mathop{\mid\kern
-1pt\joinrel\relbar\joinrel\relbar}\limits^{*}_{#1}}}
\def\longtra#1{\smash{\mathop{\mid\kern
-1pt\joinrel\relbar\joinrel\relbar\joinrel\relbar}\limits^{*}_{#1}}}
\def\vlongtra#1{\smash{\mathop{\mid\kern
-1pt\joinrel\relbar\joinrel\relbar\joinrel\relbar\joinrel\relbar}\limits^{*}_{#1}}}
\def\vvlongtra#1{\smash{\mathop{\mid\kern
-1pt\joinrel\relbar\joinrel\relbar\joinrel\relbar\joinrel\relbar\joinrel\relbar}\limits^{*}_{#1}}}
\def\vvvlongtra#1{\smash{\mathop{\mid\kern
-1pt\joinrel\relbar\joinrel\relbar\joinrel\relbar\joinrel\relbar\joinrel\relbar\joinrel\relbar}\limits^{*}_{#1}}}
\def\etra#1{\smash{\mathop{\mid\kern
-1pt\joinrel\relbar\joinrel\relbar}\limits_{#1}}}

\def\iff{\Leftrightarrow}
\def\Rw{\Rightarrow}
\def\oo{\overline}

\def\N{\mathbb{N}}

\def\at{\mbox{At}}

\def\het{\mbox{ht}}

\def\flatx{\mbox{Fl}}

\def\TR{{\cal{T}}}

\def\H{{\cal{H}}}
\def\BoR{{\cal{BR}}}

\def\p{\varphi}

\def\inv{^{-1}}

\def\bi{\begin{itemize}}
\def\ei{\end{itemize}}
\def\beq{\begin{equation}}
\def\eeq{\end{equation}}

\newtheorem{T}{Theorem}[section]
\newcommand{\bt}{\begin{T}}
\newcommand{\et}{\end{T}}
\newcommand{\ftd}{$\square$\end{T}}

\newtheorem{Proposition}[T]{Proposition}
\newcommand{\bp}{\begin{Proposition}}
\newcommand{\ep}{\end{Proposition}}
\newcommand{\fpd}{$\square$\end{Proposition}}

\newtheorem{Lemma}[T]{Lemma}
\newcommand{\bl}{\begin{Lemma}}
\newcommand{\el}{\end{Lemma}}
\newcommand{\fld}{$\square$\end{Lemma}}

\newtheorem{Corol}[T]{Corollary}
\newcommand{\bc}{\begin{Corol}}
\newcommand{\ec}{\end{Corol}}
\newcommand{\fcd}{$\square$\end{Corol}}

\newtheorem{Result}[T]{Result}
\newcommand{\br}{\begin{Result}}
\newcommand{\er}{\end{Result}}
\newcommand{\frd}{$\square$\end{Result}}

\newtheorem{Example}[T]{Example}
\newcommand{\be}{\begin{Example}}
\newcommand{\ee}{\end{Example}}

\newtheorem{Problem}[T]{Problem}
\newcommand{\bq}{\begin{Problem}}
\newcommand{\eq}{\end{Problem}}

\newcommand{\proof}
   {\par\medbreak\noindent{\bf Proof}.\enspace}

\newcommand{\qed}{
$\Box$
\par\bigbreak}

\def\abstract#1{\par\bigskip
\begingroup\small
\baselineskip=12truept
\begin{center}ABSTRACT\end{center}
\par\medskip\par\noindent
\null\hfill\hbox{\vbox{\hsize=5truein\noindent#1}}
\hfill\null\par\endgroup\par}

\title{On the lattice of flats of a boolean representable simplicial complex}
\author{{\bf Stuart Margolis, John Rhodes and Pedro V. Silva}}

\date{\today}

\begin{document}
\maketitle

\begin{center}\small
2010 Mathematics Subject Classification:  05E45, 06B15, 06B05, 05B35

\bigskip

Keywords: Simplicial complex, hereditary collection, boolean
representable, lattice of flats, atomistic lattice
\end{center}

\abstract{It is shown that the lattices of flats of boolean
  representable simplicial complexes are always atomistic, but
  semimodular if and only if the complex is a matroid. A canonical
  construction is introduced for 
  arbitrary finite atomistic lattices, providing a characterization of
the lattices of flats of boolean representable simplicial complexes
and a decidability condition. 
We remark that every finite
lattice occurs as the lattice of flats of some simplicial complex.
}

\section{Introduction}

In a series of three papers \cite{IR1,IR2,IR3}, Izhakian and Rhodes
introduced the concept of boolean representation for various algebraic
and combinatorial structures. These ideas were inspired by previous
work by Izhakian and Rowen on supertropical matrices (see
e.g. \cite{Izh2,IRow2,IRow,IRow22}), and were subsequently developed by
Rhodes and Silva in a recent monograph, devoted to boolean
representable simplicial complexes \cite{RSm}.

The original approach is to consider matrix representations over the
superboolean semiring $\mathbb{SB}$, using appropriate notions of vector
independence and rank. Writing $\mathbb{N} = \{ 0,1,2,\ldots \}$, we
can define $\mathbb{SB}$ as the quotient of $(\mathbb{N},+,\cdot)$ (usual
operations) by the congruence which identifies all integers $\geq 2$.
In this context, boolean representation refers to matrices using only
0 and 1 as entries.

In this paper, we view (finite) simplicial complexes in their abstract form, as
hereditary collections. However, boolean representations have also
provided results of a more geometric and topological nature, namely on
the homotopy of (geometric) simplicial complexes (see \cite[Chapter 7]{RSm},
\cite{MRS}). 

As an alternative, boolean representability can be characterized by
means of the lattice of flats. The lattice of flats plays a fundamental
role in matroid theory but is not usually considered for arbitrary
simplicial complexes, probably due to the fact that, unlike in the
matroid case, the structure of a simplicial complex cannot in general be
recovered from its lattice of flats. However, this is precisely what
happens with boolean representable simplicial complexes (BRSCs). If $\H\; =
(V,H)$ is a simplicial complex and $\flatx\H$ denotes its lattice of
flats, then $\H$ is boolean representable if and only if $H$ equals
the set of transversals of the successive differences for chains in
$\flatx\H$. This implies in particular that all (finite) matroids are boolean
representable. If the BRSC is simple, it can be characterized by the
isomorphism class of its lattice of flats, plus a bijection from its
vertex set onto the set of atoms of the lattice.

Therefore it is a natural question to inquire about the nature of the
lattices of flats of BRSCs. We
note that in the matroid case the lattices of flats are precisely the
geometric lattices (atomistic and semimodular) \cite[Theorem 1.7.5]{Oxl}.

We note that it is known that any finite lattice embeds in a finite
partition lattice \cite{PT}, which is a geometric lattice, so corresponds to a
matroid. Therefore every finite lattice is isomorphic to some full
sublattice of the lattice of flats of some $\H \in \BoR$.

The paper is organized as follows. In Section 2 we present all the
basic notions and results needed in the paper. In Section 3 we show
that every finite lattice occurs as the lattice of flats of some
simplicial complex. In Section 4 we show that lattices of flats of
BRSCs are always atomistic, but
not necessarily semimodular.
In Section 5 we construct a simple BRSC
$\TR_L$ for every finite atomistic lattice $L$ and show that $L \cong
\flatx\H$ for some BRSC $\H$ if
and only if $L \cong \flatx\TR_L$. Moreover, in this case $\TR_L$ is
isomorphic to some restriction of $\H$, and if $\H$ is simple, then
$\TR_L$ is isomorphic to $\H$. This provides an easy way of
deciding whether or not a finite lattice is isomorphic to the lattice
of flats of some BRSC, and the complexity is polynomial for fixed
height. In Section 
6, we prove a graph-theoretical characterization for lattices of
height 3, also decidable in polynomial time. 

\section{Preliminaries}

All lattices and simplicial complexes in this paper are assumed to be
finite. The reader is assumed to have some familiarity with basic
notions of lattice theory, being referred to \cite{Gra}.

Given a set $V$ and $n \geq 0$, we denote by $P_n(V)$
(respectively $P_{\leq n}(V)$) the set
of all subsets of $V$ with precisely (respectively at most) $n$
elements. To simplify notation, we shall often represent sets $\{ a_1, a_2,
\ldots, a_n\}$ in the form $a_1a_2\ldots a_n$.

A (finite) simplicial complex is a structure of the form $\H\; =
(V,H)$, where $V$ is a finite nonempty set and $H \subseteq 2^V$ is
nonempty and closed under taking subsets. Simplicial complexes, in
this abstract viewpoint, are also known as {\em hereditary collections}.

Two simplicial complexes $(V,H)$ and $(V',H')$ are {\em isomorphic} if
there exists a bijection $\p:V \to V'$ such that
$$X \in H \mbox{ if and only if }X\p \in H'$$
holds for every $X \subseteq V$.

If $\H\; = (V,H)$ is a simplicial complex and $W \subseteq V$ is
nonempty, we call 
$$\H|_W = (W,H \cap 2^W)$$
the {\em restriction} of $\H$ to $W$. It is obvious that $\H|_W$ is
still a simplicial complex.

We say that $X
\subseteq V$ is a {\em flat} of $\H$ if
$$\forall I \in H \cap 2^X \hspace{.2cm} \forall p \in V \setminus X
\hspace{.5cm} I \cup \{ p \} \in H.$$
The set of all flats of $\H$ is denoted by 
$\flatx\H$. 

Clearly, the intersection of any set of flats (including $V =
\cap\emptyset$) is still a flat. If we order $\flatx\H$ by inclusion,
it is then a $\wedge$-semilattice, and therefore a lattice with
$$(X \vee Y) = \cap\{ F \in \flatx\H \mid X \cup Y \subseteq F \}$$
for all $X,Y \in \flatx\H$. 
We call $\flatx\H$ the {\em lattice of flats} of $\H$.
The lattice of flats induces a closure operator on $2^V$ defined by
$$\oo{X} = \cap\{ F \in \flatx \H \mid X \subseteq F \}$$
for every $X \subseteq V$. 

We say that $X$ is a {\em transversal of the
successive differences} for a chain of subsets
$$A_0 \subset A_1 \subset \ldots \subset A_k$$
if $X$ admits an enumeration $x_1,\ldots , x_k$ such that $x_i \in A_i
\setminus A_{i-1}$ for $i = 1,\ldots,k$. 

Let $\H\; = (V,H)$ be a simplicial complex. If $X \subseteq V$ is a
{\em transversal of the successive differences} for a chain 
$$F_0 \subset F_1 \subset \ldots \subset F_k$$
in $\flatx \H$, it follows easily by induction that $\{ x_1, x_2,
\ldots x_i \} \in H$ for $i = 0,\ldots,k$. In particular, $X \in H$.

We say that $\H$ is {\em
  boolean representable} if every $X \in H$ is a transversal of the
successive differences for a chain in $\flatx \H$. We denote by $\BoR$
the class of all (finite) boolean representable simplicial
complexes (BRSCs). 

A simplicial complex $\H\; = (V,H)$ is called a {\em matroid} if it
satisfies the {\em exchange property}:
\bi
\item[(EP)]
For all $I,J \in H$ with $|I| = |J|+1$, there exists some
  $i \in I\setminus J$ such that $J \cup \{ i \} \in H$.
\ei
All matroids are boolean representable, but the converse is not true
(see \cite[Example 5.2.11]{RSm}).

\section{Arbitrary complexes}

In the literature, the lattice of flats is defined only for matroids,
but there is no reason to prevent considering it for arbitary
simplicial complexes. It turns out that every finite lattice can be
constructed this way.

\bt
\label{lsc}
Every finite lattice is isomorphic to the lattice of flats of some
simplicial complex.
\et

\proof
Let $L$ be a finite lattice. If $L$ is trivial, then $L \cong
\flatx\H$ for $\H\; = (V,\{ \emptyset\})$, hence we assume that $L$ is
nontrivial. Let $B$ and $T$ denote the bottom and the top elements of
$L$. Let
$$V = \{ a^{(i)} \mid a \in L \setminus \{ B \},\; i = 1,2,3 \}$$
and let $\pi:V \to L \setminus \{ B \}$ be the canonical mapping. Let
$$J = \{ X \subset V \; {\big{\lvert}} \; \pi|_X \mbox{ is injective}\}.$$
For every $P \subseteq L \setminus \{ B \}$, let 
$$P\alpha = \{ a \in L \setminus \{ B \} \mid p \not\leq a \mbox{ for
  every }p \in P \mbox{ and $(a \vee p) \neq q$ for all distinct }p,q
\in P\}.$$ 
We define
$$H = J \cup \{ X \cup \{ a^{(1)}, a^{(2)} \} \mid X \in J,\; a \in
X\pi\alpha \}.$$
It is easy to check that $\H\; = (V,H)$ is a (finite) simplicial
complex. This follows from $J$ being closed under taking subsets and 
$$X \cup \{ a^{(i)} \} \in J \mbox{ whenever } X \in J,\; a \in
X\pi\alpha \mbox{ and } i \in \{ 1,2\}.$$

For every $a \in L$, let 
$$a\beta = \{ x \in L \setminus \{ B \} \mid x \leq a\}.$$
Since $P_1(V) \subseteq H$, we get 
\beq
\label{lsc1}
\oo{\emptyset} = \emptyset = B\beta\pi\inv.
\eeq
Next let $a \in L \setminus \{ B \}$ and $i \in \{ 1,2,3\}$. We show
that
\beq
\label{lsc2}
\oo{a^{(i)}} = a\beta\pi\inv.
\eeq

Since $a^{(i)} \in H$ and $\{ a^{(1)},a^{(3)} \}, \{ a^{(2)},a^{(3)}
\} \notin H$, we get $a^{(3)} \in \oo{a^{(i)}}$ and thus $a\pi\inv
\subseteq \oo{a^{(i)}}$. Now let $B \neq x < a$ in $L$ and let $j \in \{
1,2,3\}$. Then $\{ a^{(1)},a^{(2)}, x^{(j)} \} \notin H$ because $a
\notin x\alpha$. Since $\{ a^{(1)},a^{(2)} \} \in H$, we get $x^{(j)}
\in \oo{a^{(i)}}$ and so $a\beta\pi\inv \subseteq \oo{a^{(i)}}$. 

Thus it is enough to prove that $a\beta\pi\inv \in \flatx\H$. Let $I
\subseteq a\beta\pi\inv$ be such that $I \in H$ and let $y^{(j)} \in V
\setminus a\beta\pi\inv$. If $I \in J$, then $y^{(j)} \notin
a\beta\pi\inv$ together with $I \subseteq a\beta\pi\inv$ yield $I
\cup \{ y^{(j)} \} \in J \subseteq H$, hence we may assume that $I =
X \cup \{ a^{(1)}, a^{(2)} \}$ for some $X \in J$ and $a \in
X\pi\alpha$. 

Suppose that $x^{(k)} \in X$. Then $x^{(k)} \in I \subseteq
a\beta\pi\inv$ and so $x \leq a$. On the other hand, $a \in
X\pi\alpha$ yields $x \not\leq a$, a contradiction. Thus $X =
\emptyset$ and so $I = \{ a^{(1)}, a^{(2)} \}$. 
It follows that
$$I \cup \{ y^{(j)} \} = \{ y^{(j)} \} \cup \{ a^{(1)},
a^{(2)} \}.$$
Clearly, $\{ y^{(j)} \} \in J$ and $y^{(j)} \notin a\beta\pi\inv$
yields $a \in y^{(j)}\pi\alpha$. Therefore $I \cup \{ y^{(j)} \} \in
H$. It follows that $a\beta\pi\inv \in \flatx\H$ and so (\ref{lsc2})
holds.

Next we show that
\beq
\label{lsc3}
\oo{\{ a^{(i)},b^{(j)} \} } = (a\vee b)\beta\pi\inv
\eeq
holds for all $a,b \in L \setminus \{ B \}$.

Let $c = (a\vee b)$. In view of (\ref{lsc2}), we may assume that $c >
a,b$. It is easy to check that $\{ a^{(1)}, a^{(2)}, b^{(1)} \} \in H$
(because $b \not\leq a$)
and $\{ a^{(1)}, a^{(2)}, b^{(1)}, c^{(1)} \} \notin H$ (because $c =
(a\vee b)$). Since $\{ a^{(1)}, a^{(2)}, b^{(1)} \} \subseteq \oo{\{
  a^{(i)},b^{(j)} \} }$ by (\ref{lsc2}), we get $c^{(1)} \in \oo{\{
  a^{(i)},b^{(j)} \} }$ and so $c\beta\pi\inv \subseteq \oo{\{
  a^{(i)},b^{(j)} \} }$ by (\ref{lsc2}). Now the opposite
inclusion follows from $\{ a^{(i)},b^{(j)} \} \subseteq (a\vee
b)\beta\pi\inv$, therefore (\ref{lsc3}) holds.

Now we claim that
\beq
\label{lsc4}
\flatx\H \; = \{ a\beta\pi\inv \mid a \in L\}.
\eeq
The opposite inclusion follows from (\ref{lsc1}) and
(\ref{lsc2}). Take $F \in \flatx\H$. In view of (\ref{lsc1}), we may
assume that $F \neq \emptyset$. It follows from (\ref{lsc3}) that
$F\pi$ has a maximum $a \neq B$, and (\ref{lsc2}) yields $F = a\beta\pi\inv$.
Thus (\ref{lsc4}) holds and it follows easily that
$$\begin{array}{rcl}
L&\to&\flatx\H\\
a&\mapsto&a\beta\pi\inv
\end{array}$$
is an isomorphism of posets and therefore a lattice isomorphism.
\qed

\section{A necessary condition}

In this section, we start the discussion on the lattices of flats of BRSCs.

The main theorem of the section proves a necessary condition for a lattice to be
isomorphic to such a lattice of flats. This implies
that we cannot assume the simplicial complexes to be
boolean representable in Theorem \ref{lsc}. But first we prove 
two simple lemmas.

\bl
\label{flatrest}
Let $\H\; = (V,H)$ be a simplicial complex and let $W \subseteq
V$. Let $F \in {\rm Fl}\H$. Then $F \cap W \in {\rm Fl}(\H|_W)$.
\el

\proof
Let $H' = H \cap 2^W$. Assume that $I \in H' \cap 2^{F\cap W}$ and $p
\in W \setminus (F\cap W)$. Then $I \in H \cap 2^{F}$ and $p
\in V \setminus F$, hence $F \in \flatx\H$ yields $I \cup \{ p \} \in
H$. Since $I \cup \{ p \} \subseteq W$, we get $I \cup \{ p \} \in
H'$. Thus $F\cap W \in \flatx(\H|_W)$.
\qed

\bl
\label{pai}
Let $\H \; = (V,H)$ be a simplicial complex and let 
$V' = \{ p \in V \mid \{ p \} \in H\}$. Then ${\rm Fl}\H \cong
{\rm Fl}(\H|_{V'})$.
\el

\proof
Write $\H' = \H|_{V'}$ and $H' = H \cap 2^{V'}$. 
Let $\alpha:\flatx\H \to \flatx\H'$ and $\beta:\flatx\H' \to \flatx\H$
be the mappings defined by
$$F\alpha = F \cap V',\quad F'\beta = F' \cup (V \setminus V').$$
By Lemma \ref{flatrest}, $\alpha$ is well defined. Now let $F' \in
\flatx\H'$. Assume that $I \in H \cap 2^{F'\beta}$ and $p 
\in V \setminus F'\beta$. By definition of $V'$, we have $I \subseteq
V'$, hence $I \in H' \cap 2^{F'}$. Since also $p
\in V' \setminus F'$, then $F' \in \flatx\H'$ yields $I \cup \{ p \} \in
H' \subseteq H$. Thus $F'\beta \in \flatx\H$ and $\beta$ is well defined.

Next we show that the mappings $\alpha$ and $\beta$ are mutually
inverse.
Let $F \in \flatx\H$. Then 
$$F\alpha\beta = (F \cap V') \cup (V \setminus V')$$
and so $F \subseteq F\alpha\beta$. Suppose that $p \in F\alpha\beta
\setminus F$. Then $p \in V \setminus V'$. Since $\emptyset \in H \cap
2^F$ and $\{ p \} \notin H$, then $F \in \flatx\H$ yields $p \in
F$, a contradiction. Thus $F\alpha\beta = F$ and so $\alpha\beta = 1$.

On the other hand, for every $F' \in \flatx\H'$ we have
$$F'\beta\alpha =  (F' \cup (V \setminus V'))\cap V' = F',$$
thus $\beta\alpha = 1$. Therefore $\alpha$ and $\beta$ are mutually
inverse.

Since $\alpha$ and $\beta$ are both order-preserving, they are poset
isomorphisms and therefore lattice isomorphisms.
\qed



Let $L$ be a lattice. Given $a,b \in L$, we say that $b$ {\em covers}
$a$ if $a < b$ and there is no $c \in L$ such that $a < b < c$.
An {\em atom} of $L$ is a element covering the bottom
element $B$. We denote by $\at(L)$ the set of atoms of $L$.
The lattice $L$ is {\em atomistic} if every element of
$L$ is a join of atoms. 
We show next that being atomistic is a necessary condition for a
lattice to be isomorphic to the lattice of flats of a BRSC.

\bt
\label{latomi}
Let $\H\; = (V,H) \in \BoR$. Then:
\bi
\item[(i)] if $P_1(V) \subseteq H$, then $V$ is the union of the atoms
  of ${\rm Fl}\H$;
\item[(ii)] ${\rm Fl}\H$ is atomistic.
\ei
\et

\proof
(i) In \cite{MRS}, we introduced the equivalence relation on $V$
defined by
$$p\eta q \quad \mbox{ if }\oo{p} = \oo{q}$$
and the {\em simplification} $\H_S =
(V/\eta,H/\eta)$. By \cite[Proposition 4.2]{MRS}(iii), we have $\flatx\H
\cong \flatx\H_S$. On the other hand, if $\p:V \to V/\eta$ denotes the
canonical projection, it follows from
\cite[Proposition 4.2]{MRS}(ii) that the atoms of $\flatx\H$ are of
the form $A\p\inv$, where $A$ is an atom of $\H_S$. But the the atoms
of $\flatx\H_S$ are the singleton sets, thus we are done.

(ii) By Lemma \ref{pai}, we may assume that $P_1(V) \subseteq H$, and
by \cite[Proposition 4.2]{MRS}(iii) we may assume that $\H$ is
simple. Hence the atoms of $\H$ are the singleton sets $\{ p \}$ with
$p \in V$. Now the claim becomes obvious. Note that $\emptyset$ is the
bottom element of $\flatx\H$, which is 
the join of the empty set of atoms.
\qed

\bc
\label{smallest}
The 3-point chain
$$\xymatrix{
\bullet \ar@{-}[d] \\
\bullet \ar@{-}[d] \\
\bullet
}$$
is the smallest lattice not isomorphic to the lattice of flats of some
$\H\; \in \BoR$. 
\ec

\proof
In fact, this chain happens to be the smallest non atomistic
lattice. On the other hand, the trivial
lattice and the 2-point lattice occur as the lattices of flats for
$(V, \{ \emptyset \})$ and $(V,P_{\leq 1}(V))$, respectively. 
Now the claim follows from Theorem \ref{latomi}(ii).
\qed

A lattice $L$ is said to be (upper)  
{\em semimodular} if has no sublattice of the form 
$$
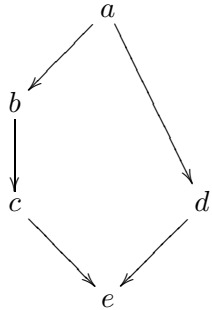
\xymatrix{
& a \ar[dl] \ar[ddr] & \\
b \ar[d] && \\
c \ar[dr] && d \ar[dl] \\
& e &
}$$
with $d$ covering $e$. A {\em geometric} lattice is a lattice which is
both semimodular and atomistic. It is well known that a finite lattice
is geometric if and only if it is isomorphic
to the lattice of flats of some (finite) matroid \cite[Theorem 
 1.7.5]{Oxl}.

The next example shows that the lattice of flats of a BRSC is not
necessarily semimodular (we must take a BRSC which is not a matroid). 

\be
\label{exnsm}
Let $V = \{ 1,2,3,4\}$ and $H = P_{\leq 2}(V) \cup \{ 123, 124 \}$. 
Then $\H\; = (V,H) \in \BoR$ but $\flatx\H$ is not semimodular. 
\ee

Indeed, it is easy to check that $\flatx\H\; = P_{\leq 1}(V) \cup \{
V, 12 \}$, and can therefore be described as
$$
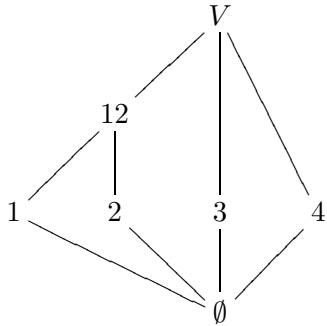
\xymatrix{
&& V \ar@{-}[dl] \ar@{-}[dd] \ar@{-}[ddr] & \\
&12 \ar@{-}[dl] \ar@{-}[d] && \\
1 \ar@{-}[drr] & 2 \ar@{-}[dr] & 3 \ar@{-}[d] & 4 \ar@{-}[dl] \\
&& \emptyset &
}$$
and is therefore boolean representable (see \cite[Example 5.2.11]{RSm}).

By considering the sublattice $\{ V, 12, 2, 4,
\emptyset \}$, we deduce that $\flatx\H$ is not semimodular.

\section{The minimal simplicial complex on a lattice of flats}

We show in this section that, whenever a lattice is isomorphic to the lattice
of flats of some $\H\; \in \BoR$, there exists a minimal simplicial
complex satisfying this condition. To prove this 
claim, we introduce the following construction. 

Let $L$ be an atomistic lattice. 
For every $x \in L$, let 
$$x\xi = \{ a \in \at(L) \mid a \leq x\}.$$
We say that $A \subseteq \at(L)$ is a {\em transversal of a chain} in $L$
if there exists an enumeration $a_1, \ldots, a_m$ of the elements of
$A$ and a chain
$$x_0 < x_1 < \ldots < x_m$$
in $L$ such that $a_i \in x_i\xi \setminus x_{i-1}\xi$ for $i =
1,\ldots,m$. 

Let $T_L \subseteq 2^{\at(L)}$ consist of all the transversals of some
chain in $L$ and write  
$\TR_L = (\at(L),T_L)$. It is immediate that $\TR_L$ is a simplicial
complex. We prove the following lemma.

\bl
\label{atbr}
Let $L$ be an atomistic lattice. Then:
\bi
\item[(i)] $x\xi \in {\rm Fl}\TR_L$ for every $x \in L$;
\item[(ii)] $\xi$ is a lattice isomorphism of $L$ onto $L\xi \subseteq
{\rm Fl}\TR_L$;
\item[(iii)] $\TR_L \in \BoR$;
\item[(iv)] $\TR_L$ is simple. 
\ei
\el

\proof
(i) Let $x \in L$. If $x = T$, then the conclusion is clear, so assume
that $x \neq T$. 
Assume that $A \in T_L \cap 2^{x\xi}$ and $p \in \at(L)
\setminus x\xi$. Since $A \in T_L$, there exists an enumeration $a_1,
\ldots, a_m$ of $A$ and a chain
$x_0 < \ldots < x_m$
in $L$ such that $a_i \in x_i\xi \setminus x_{i-1}\xi$ for $i =
1,\ldots,m$. Now
\beq
\label{emb2}
(x_0 \wedge x) \leq (x_1 \wedge x)  \leq \ldots \leq (x_m \wedge
x).
\eeq
Since $A \subseteq x\xi$, we get
$a_i \in (x_i \wedge x)\xi \setminus (x_{i-1} \wedge x)\xi$ for $i =
1,\ldots,m$, hence the ordering in (\ref{emb2}) must be strict. Since 
$(x_m \wedge
x)\xi \subseteq x\xi$, then $p \in T\xi \setminus (x_m \wedge
x)\xi$ and so $A \cup \{ p \}$ is a transversal of the chain
$$(x_0 \wedge x) < (x_1 \wedge x) < \ldots < (x_m \wedge
x) < T$$
in $K$. Thus $A \cup \{ p \} \in T_L$ and so $x\xi \in \flatx\TR_L$.

(ii) Since $L$ is atomistic, we have 
\beq
\label{derby2}
x = \vee (x\xi) \mbox{ for every }x \in L,
\eeq
hence $\xi$ is injective. 

Let $x,y \in L$. It is immediate that $x \leq y$ implies $x\xi \subseteq
y\xi$, hence $\xi$ is order-preserving. Finally, in view of
(\ref{derby2}), $x\xi \subseteq y\xi$ yields 
$$x = \vee(x\xi) \leq \vee(y\xi) = y,$$
hence $\xi:L \to L\xi$ is a poset isomorphism and therefore a lattice
isomorphism (preserving top and bottom).  

(iii) Let $A
\in T_L$. Then there exists an enumeration $a_1, \ldots, a_m$ of the elements of
$A$ and a chain
$$x_0 < x_1 < \ldots < x_m$$
in $L$ such that $a_i \in x_i\xi \setminus x_{i-1}\xi$ for $i =
1,\ldots,m$. Since $x_i\xi \in \flatx\TR_L$ for $i = 0,\ldots,m$, it
follows that $A$ is a transversal of the
successive differences for a chain in $\flatx \TR_L$.
Therefore $\TR_L$ is boolean representable.

(iv) It is immediate that $P_2({\rm At}(L)) \subseteq T_L$.
\qed

The following theorem asserts the canonical role played by $\TR_L$.

\bt
\label{canonic}
Let $L$ be an atomistic lattice. Then the following conditions are
equivalent.
\bi
\item[(i)] $L \cong {\rm Fl}\H$ for some $\H \in \BoR$;
\item[(ii)] $L \cong {\rm Fl}\TR_L$;
\item[(iii)] $\xi:L \to {\rm Fl}\TR_L$ is onto;
\item[(iv)] $(\vee F)\xi \subseteq F$ for every $F \in {\rm Fl}\TR_L$.
\ei
Moreover, in this case $\TR_L$ is isomorphic to some restriction of
$\H$, and for $\H$ simple we get an isomorphism.
\et

\proof
(i) $\Rw$ (ii). Let $\H\; = (V,H)$.
By Lemma \ref{pai}, we may assume that $P_1(V)
\subseteq H$, and by \cite[Proposition 4.2]{MRS}(iii) we may assume
that $\H$ is simple (replacing $\H$ by its simplification $\H_S$,
isomorphic to some restriction of $\H$). 

Let $\p:L \to \flatx\H$ be a lattice isomorphism. Since $\at(\flatx\H)
= \{ \{ p \} \mid p \in V\}$, $\p$ induces a bijection $\p':\at(L)
\to V$ defined by $a\p = \{ a\p'\}$.
We claim that
\beq
\label{canonic1}
\TR_L \cong \H.
\eeq

Indeed, let $X \subseteq \at(L)$. Then $X \in T_L$ if and only if
there exists an enumeration $a_1,\ldots,a_k$ of the elements of $X$
and a chain $\ell_0 < \ell_1 < \ldots < \ell_k$ in $L$ such that $a_i
\in \ell_i\xi \setminus \ell_{i-1}\xi$ for $i = 1,\ldots,k$. It is
easy to check that $\ell_0\p < \ell_1\p < \ldots < \ell_k\p$ is a chain in
$\flatx\H$ such that $a_i\p' \in \ell_i\p \setminus \ell_{i-1}\p$ for
$i = 1,\ldots,k$. Hence $X\p' \in H$. Similarly, $X\p' \in H$ implies
$X \in T_L$, therefore $\p'$ defines an isomorphism between $\TR_L$
and $\H$, and so (\ref{canonic1}) holds.

Note that, for arbitrary $\H$, the reductions performed above through 
Lemma \ref{pai} and \cite[Proposition 4.2]{MRS}(iii) replace the
original BRSC $\H$ by one of its restrictions, so $\TR_L$ is
indeed isomorphic to some restriction of $\H$.

(ii) $\Rw$ (i). By Lemma \ref{atbr}(i).

(ii) $\iff$ (iii). In view of Lemma \ref{atbr}(ii), and since $L$ is
finite.

(iii) $\Rw$ (iv). Let $F \in \flatx\TR_L$. By condition (iii), we have
$F = x\xi$ 
for some $x \in L$. Hence $a \leq x$ for every $a \in F$ and so $\vee
F \leq x$. Thus $(\vee F)\xi \subseteq x\xi = F$.

(iv) $\Rw$ (iii). Let $F \in \flatx\TR_L$ and let $x = \vee F$. Since
$a \in F$ implies $a \leq x$, we have $F \subseteq x\xi$. On the other
hand, condition (iv) yields $x\xi \subseteq F$, hence $\xi:L \to
\flatx\TR_L$ is onto.
\qed

\bc
\label{deci}
Let $L$ be a lattice. Then it is decidable whether or not
$L \cong {\rm Fl}\H$ for some $\H\; \in \BoR$.
\ec

\proof
By Theorem \ref{latomi}(ii), $L$ being atomistic is a necessary
condition. Since being atomistic is certainly decidable, we may assume
that $L$ is atomistic. 

Since we can successively compute $\TR_L$ and $\flatx(\TR_L)$, the
claim now follows from Theorem \ref{canonic}.
\qed

Predictably, being atomistic is not a sufficient condition, as
we show in  
the next result. We recall that the {\em height} of a lattice $L$,
denoted by $\het(L)$, is the maximal length $n$ of a chain
$x_0 < x_1 < \ldots < x_n$ in $L$. Given a set $X$, let
$(2^X,\subseteq)$ be the lattice of all subsets of $X$ with respect to
inclusion.

\bp
\label{athe}
Let $L$ be a lattice with ${\rm ht}(L) = |{\rm At}(L)|$. Then the
following conditions are equivalent:
\bi
\item[(i)] $L \cong {\rm Fl}\H$ for some $\H\; \in \BoR$;
\item[(ii)] $L \cong (2^{{\rm At}(L)},\subseteq)$.
\ei
\ep

\proof
(i) $\Rw$ (ii).
By Theorem \ref{latomi}(ii), $L$ is atomistic. Let
$$x_0 < x_1 <
\ldots < x_n$$ 
be a chain of maximal length in $L$. Since $L$ is atomistic, for every $i \in \{
1,\ldots,n\}$ there exists some $a_i \in x_i\xi \setminus
x_{i-1}\xi$. It is easy to see that the $a_i$ are all distinct and so
$\het(L) = |\at(L)|$ yields $\at(L) = \{ a_1, \ldots,a_n\}$. Hence
$\at(L) \in T_L$ and so $T_L = 2^{\at(L)}$. Thus $2^{\at(L)} =
\flatx\TR_L$ is a lattice (with height $|\at(L)|$) with respect to 
inclusion. Now the claim follows from Theorem \ref{canonic}.

(ii) $\Rw$ (i). Let $\H\; = (\at(L),2^{\at(L)})$. Then $\flatx\H\; =
2^{\at(L)}$ and so $\H\in\BoR$. Since $L \cong \flatx\H$, we
are done. 
\qed

We can now produce the equivalent of Corollary \ref{smallest} for
atomistic lattices.

\bc
\label{smallestat}
The lattice $L$ depicted by
$$
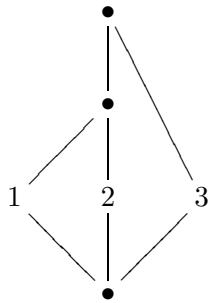
\xymatrix{
& \bullet \ar@{-}[d] \ar@{-}[ddr] & \\
& \bullet \ar@{-}[dl] \ar@{-}[d] & \\
1 \ar@{-}[dr] & 2 \ar@{-}[d] & 3 \ar@{-}[dl] \\
& \bullet &
}$$
is the smallest atomistic lattice not isomorphic to the lattice of flats of some
$\H\; \in \BoR$. 
\ec

\proof
Since the above lattice, which is clearly atomistic, has 3 atoms,
height 3 and 6 elements, it is not isomorphic to the lattice of flats
of some $\H\in \BoR$ by Proposition \ref{athe}.

Indeed, it is easy to check that $\TR_L = (123, P_{\leq 3}(123))$ and
so $\flatx\TR_L$ is the lattice
$$
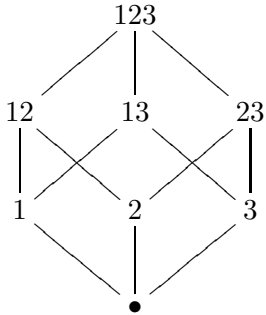
\xymatrix{
& 123 \ar@{-}[d] \ar@{-}[dr] \ar@{-}[dl] & \\
12 \ar@{-}[d] \ar@{-}[dr] & 13 \ar@{-}[dl] \ar@{-}[dr] 
& 23 \ar@{-}[d] \ar@{-}[dl]\\
1 \ar@{-}[dr] & 2 \ar@{-}[d] & 3 \ar@{-}[dl] \\
& \bullet &
}$$
clearly not isomorphic to $L$.

It is easy to check that any other atomistic lattice with less than
7 elements must have height $\leq 2$ and is therefore geometric. It
follows that it must be isomorphic to the lattice of flats of a
matroid, and matroids are boolean representable.
\qed

We end this section by showing that the complexity of the algorithm
outlined in Corollary \ref{deci} is polynomial for
fixed height.

We recall the $O$ notation from complexity theory. Let $P$ be an
algorithm defined for instances depending on parameters
$n_1,\ldots,n_k$. If $\p:\N^k \to \N$ is a function, we write $P \in
O((n_1,\ldots,n_k)\p)$ if there exist constants $K,L > 0$ such that
$P$ processes each instance of type $(n_1,\ldots,n_k)$ in time $\leq
K((n_1,\ldots,n_k)\p) + L$ (where time is measured as the number of
elementary operations performed).

\bp
\label{polc}
It is decidable in time $O(n^{3h})$ whether an atomistic lattice of
height $h$ with $n$ atoms is isomorphic to ${\rm Fl}\H$ for some $\H\;
\in \BoR$. 
\ep

\proof
Let $M = (m_{\ell a})$ be the $L \times \at(L)$ defined by
$$m_{\ell a} = \left\{
\begin{array}{ll}
0&\mbox{ if } x \geq a\\
1&\mbox{ otherwise}
\end{array}
\right.$$
It follows from the results in \cite[Section 3.5]{RSm} that $M$ is a
boolean matrix representation of $\TR_L$. 

We may assume $h \geq 3$. By \cite[Theorem 7.4]{MRS}, it is possible
to compute in time $O(n^{3d+3})$ the 
list of flats of a simplicial complex of dimension $d$ defined by
a reduced boolean matrix with $n$ columns. Since $L$ has height $h$,
the dimension of $\TR_L$ (corresponding to the maximum size of an
element of $T_L$ minus 1) is $h-1$. Moreover, $M$ is reduced (all rows
are distinct) since $L$ is atomistic. Thus we can enumerate in time
$O(n^{3h})$ the list of flats of $\TR_L$ and compute $|\flatx\TR_L|$.

Now, in view of Lemma \ref{atbr}(ii) and Theorem \ref{canonic}, we
have $L \cong {\rm Fl}\H$ for some $\H\;
\in \BoR$ if and only if $|\flatx\TR_L| = |L|$. Thus we obtain the
claimed complexity.
\qed

\section{Height 3}

If $L$ is an atomistic lattice of height 2, then $L \cong \flatx\TR_L$ 
since $T_L = P_{\leq 2}(\at(L))$ and so $\flatx\TR_L = P_{\leq 1}(\at(L)) 
\cup \{ \at(L)\}$.
 
We provide in this section a necessary and sufficient
graph-theoretical condition for the case of (atomistic) lattices of
height 3 to be lattices of flats of a BRSC. Note that lattices of
height 3 are important since every 
simplicial complex of dimension 2 has a lattice of flats of height
$\leq 3$, and dimension 2 is a broad universe. For instance,
any finitely presented group can occur as
the fundamental group of a simplicial complex of dimension 2 (see
\cite[Theorem 7.45]{Rot}).

Let $\Gamma = (V,E)$ be a finite (undirected) graph. A {\em
clique} of $\Gamma$ is a subset $W$ of $V$ inducing a complete
subgraph of $\Gamma$. The clique $W$ is nontrivial if $|W|
\geq 2$. A nontrivial clique $W$ is a {\em superclique} if,
for all $a,b \in W$ distinct, every $c \in V \setminus W$ is not adjacent to
either $a$ or $b$. In particular, every superclique
is a maximal clique.

Given an atomistic lattice $L$ of height 3 (with top element $T$), we
define a graph $\Gamma_L = (\at(L),E_L)$ by
$$E_L = \{ \{ a,b \} \mid a,b \in \at(L),\; (a\vee b) = T \}.$$
We remark that, if $L$ is the lattice of flats of a BRSC $\H$, then $\Gamma_L$
ia actually the complement graph of $\Gamma\flatx\H$, the {\em graph
  of flats} of $\H$, introduced in \cite[Section 6.4]{RSm}.

\bt
\label{heth}
Let $L$ be a lattice of height 3. Then the following conditions are
equivalent.
\bi
\item[(i)] $L \cong {\rm Fl}\H$ for some $\H \in \BoR$;
\item[(ii)] $L$ is atomistic and $\Gamma_L$ has no supercliques.
\ei
\et

\proof
(i) $\Rw$ (ii). By Theorem \ref{latomi}, $L$ is
atomistic. By Theorem \ref{canonic}, we may assume that $\H\; = \TR_L
= (\at(L),T_L)$. 

Suppose that $W \subseteq \at(L)$ is a superclique of
$\Gamma_L$. Since $L$ has height 3, we have $W \subset \at(L)$. We
claim that $W \in \flatx\TR_L$. Indeed, let $I \in T_L \cap 2^W$ be
nonempty and
$p \in \at(L) \setminus W$. Then $I$ admits
an enumeration $a_1, \ldots, a_m$ such that
$$a_1 < (a_1 \vee a_2) < \ldots < (a_1 \vee \ldots \vee a_m).$$
Since $a_1 \edge a_2$ is an edge of $\Gamma_L$, we get $(a_1 \vee a_2)
= T$ and so $|I| \leq 2$. By Lemma \ref{atbr}(iv), we may
assume that $|I| = 2$. Since $W$ is a superclique of
$\Gamma_L$, there exists some $a \in W$ such that $\{ a,p \} \notin
E_L$. Hence $(a \vee p) < T$. Writing $I = \{ a,b \}$, we have $(a
\vee b) = T$, hence $a < (a \vee p) < T = (a \vee b \vee p)$ and so $I
\cup \{ p \} \in T_L$. Thus $W \in \flatx\TR_L$.

It follows from Theorem \ref{canonic} that $\at(L) = T\xi = (\vee
W)\xi \subseteq W$, contradicting $W \subset \at(L)$. Therefore
$\Gamma_L$ has no supercliques. 

(ii) $\Rw$ (i). Suppose that (i) fails. By Theorem \ref{canonic},
there exists some $F \in {\rm Fl}\TR_L$ such that
$(\vee F)\xi \not\subseteq F$. It follows that $|F| \geq 2$. 

Suppose that $(a \vee b) < T$ for some distinct $a,b \in F$. Let
$p \in (a \vee b)\xi \setminus F$. By Lemma \ref{atbr}(iv), we
get $abp \in T_L$. Since $(a \vee b \vee p) = (a \vee b) < T$, this
contradicts $\het(L) = 3$. Thus $(a \vee b) = T$ for all distinct $a,b
\in F$ and so $F$ is a clique of $\Gamma_L$.

Let $a,b \in F$ be distinct and let $p \in \at(L) \setminus F$. By
Lemma \ref{atbr}(iv), we get $abp \in T_L$. Now $(a \vee b) =
T$ implies that either $(a \vee p) < T$ or $(b \vee p) < T$, hence $\{
\{a, p\}, \{ b,p\} \} \not\subseteq E_L$ and so $F$ is a superclique of
$\Gamma_L$, a contradiction. Therefore (ii) holds as required.
\qed

We have remarked before that the lattices of flats of matroids are precisely the
geometric lattices. Since every matroid is boolean representable, it
follows that $\Gamma_L$ has no supercliques when $L$ is a geometric
lattice of height 3. Indeed, in this case the graph $\Gamma_L$ has no
edges at all. Suppose that $a,b \in at(L)$  satisfy $(a \vee b) =
T$. Since a geometric lattice must satisfy the Jordan-Dedekind
condition (all maximal chains have the same length), there exists some
$c \in L$ such that $a < c < T$. But then 
$$
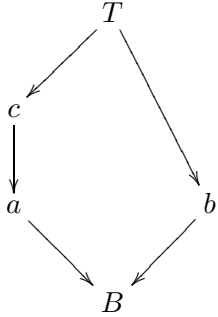
\xymatrix{
& T \ar[dl] \ar[ddr] & \\
c \ar[d] && \\
a \ar[dr] && b \ar[dl] \\
& B &
}$$
is a subsemilattice of $L$, contradicting semimodularity.

Note also that, for the lattice $L$ of Corollary \ref{smallestat},
$\Gamma_L$ is the graph
$$1 \edge 3 \edge 2$$
and has therefore supercliques ($13$ and $23$). Therefore Theorem
\ref{heth} provides an alternative way of showing that $L$ is not
isomorphic to the lattice of flats of some $\H\; \in \BoR$. 

We can show that the algorithm implicit in Theorem \ref{heth} has
polynomial complexity. 

\bp
\label{3com}
It is decidable in time $O(n^5)$ whether a graph with $n$ vertices
has no supercliques.
\ep

\proof
Let $\Gamma = (V,E)$ be a graph with $n$ vertices. Given an edge $a
\edge b$ in $\Gamma$, we define a sequence of sets of vertices
\beq
\label{3com2}
ab = W_0 \subset W_1
\subset \ldots \subset W_k \subseteq V
\eeq
as follows. 

Assume that $W_i$ is defined. If there exists some $v \in V \setminus
W_i$ adjacent to at least two vertices in $W_i$, let $W_{i+1} = W_i
\cup \{ v \}$. Otherwise, the sequence terminates at $W_i$. 

Note that this procedure is nondeterministic, but it turns out to be
confluent. Indeed, let $ab = W'_0 \subset \ldots \subset W'_{\ell}$ be
an alternative sequence. Write 
$W'_j\setminus W'_{j-1} = \{ w'_j\}$. Let $j \in \{ 1,\ldots,\ell\}$
and assume that $\{ w'_1,\ldots,w'_{j-1} \} \subseteq W_k$. Then
$w'_j$ is adjacent to at least two vertices of 
$W_k$. Since (\ref{3com2}) terminates at $W_k$, it follows that $w'_j
\in W_k$. By induction, we get $W'_{\ell} \subseteq W_k$ and so
$W'_{\ell} = W_k$ by symmetry. Therefore the procedure is confluent
and we may define $C(ab) = W_k$. 

We show that 
\beq
\label{3com1}
\mbox{every superclique of $\Gamma$ is of the form $C(ab)$ for some
  $ab \in E$.}
\eeq

Suppose that $W \subseteq V$ is a superclique of $\Gamma$. Since $|W|
\geq 2$, there exists an edge $ab \in E \cap P_2(W)$. We show that $W
= C(ab)$. 

Consider the sequence (\ref{3com2}) with $W_k = C(ab)$. Straightforward
induction shows that $W_i \subseteq W$ for $i = 0,\ldots,k$. Thus
$C(ab) \subseteq W$ and so $C(ab)$ is a clique of $\Gamma$, actually a
superclique since every $v \in V \setminus C(ab)$ is adjacent to at
most one vertex of $C(ab)$. Since every superclique is a maximal
clique, it follows that $W = C(ab)$ and so (\ref{3com1})
holds.

It is easy to see that 
\beq
\label{3com3}
\mbox{$C(ab)$ is a superclique of $\Gamma$ if and only if it is a clique.}
\eeq

Now $\Gamma$ has $O(n^2)$ edges. For each one of these edges, say
$ab$, we can compute $C(ab)$ in time $O(n^3)$. Indeed, $k \leq n-2$,
and we can go from $W_{i-1}$ to $W_i$ in time $O(n^2)$. This can be
argued by considering the adjacency matrix of $\Gamma$ (an $n \times
n$ boolean matrix): assuming that the rows and columns corresponding
to $W_{i-1}$ are marked, it 
suffices to go through the entries of the matrix once to search for the
vertex of $W_i \setminus W_{i-1}$. Using the adjacency matrix, we can
also check if $C(ab)$ is a clique in time $O(n^2)$.

In view of (\ref{3com1}) and (\ref{3com3}), we can compute in time
$O(n^5)$ any possible supercliques of $\Gamma$.
\qed

This complexity bound can most probably be improved.

It is clear that the construction of $\Gamma_L$ from $L$ can be
performed in time at most $O(n^5)$, so we get the polynomial
complexity claim on the lattice (in $O(n^5)$). Note that the general
algorithm from Proposition \ref{polc} yields only a complexity in $O(n^9)$.



The following question follows naturally from the results in this
paper: is there a more
efficient algorithm to decide if a lattice is the lattice of flats of
a BRSC?

\section*{Acknowledgements}

The first author was partially supported by Binational Science
Foundation grant number 2012080.

The second author thanks the Simons Foundation-Collaboration Grants
for Mathematicians for travel grant $\#$313548.

The third author was partially supported by
CNPq (Brazil) through a BJT-A grant (process 313768/2013-7) and
CMUP (UID/MAT/00144/2013),
which is funded by FCT (Portugal) with national (MEC) and European
structural funds through the programs FEDER, under the partnership
agreement PT2020.

\bigskip

{\sc Stuart Margolis, Department of Mathematics, Bar Ilan University,
  52900 Ramat Gan, Israel} 

{\em E-mail address:} margolis@math.biu.ac.il

\bigskip

{\sc John Rhodes, Department of Mathematics, University of California,
  Berkeley, California 94720, U.S.A.}

{\em E-mail addresses}: rhodes@math.berkeley.edu, BlvdBastille@aol.com

\bigskip

{\sc Pedro V. Silva, Centro de
Matem\'{a}tica, Faculdade de Ci\^{e}ncias, Universidade do
Porto, R. Campo Alegre 687, 4169-007 Porto, Portugal}

{\em E-mail address}: pvsilva@fc.up.pt

\end{document}